\documentclass[12pt,reqno]{amsart}

\usepackage{amsthm}
\usepackage{amsmath}
\usepackage{amsfonts}
\usepackage{amssymb}
\usepackage{hyperref}
\usepackage{fullpage}

\newtheorem{theorem}{Theorem}
\newtheorem{lemma}[theorem]{Lemma}

\theoremstyle{definition}

\theoremstyle{remark}

\numberwithin{equation}{section}

\renewcommand{\b}{\mathfrak{b}}
\renewcommand{\phi}{\varphi}
\newcommand{\h}{\mathfrak{h}}
\newcommand{\g}{\mathfrak{g}}
\newcommand{\gl}{\mathfrak{gl}}
\newcommand{\C}{\mathbb{C}}

\newcommand{\alt}{\operatorname{Alt}_s}
\newcommand{\End}{\operatorname{End}}

\allowdisplaybreaks

\title[Quantization of zero-weight super dynamical $r$-matrices]{On the quantization of zero-weight super dynamical $r$-matrices}

\author{Gizem Karaali}
\address{Department of Mathematics, Pomona College, Claremont, CA
91711}
\email{gizem.karaali@pomona.edu}

\begin{document}

\begin{abstract}
Solutions of the classical dynamical Yang-Baxter equation on a Lie superalgebra are called \emph{super dynamical $r$-matrices}. A super dynamical $r$-matrix $r$ satisfies the \emph{zero weight condition} if:
\begin{equation*}
[h\otimes 1 + 1 \otimes h, r(\lambda)] = 0 \textmd{ for all } h \in \h, \lambda \in \h^*.
\end{equation*}
\noindent
In this note we explicitly quantize zero-weight super dynamical $r$-matrices with zero coupling constant. We also answer some questions about super dynamical $R$-matrices. In particular we offer some support for one particular interpretation of the super Hecke condition. 
\end{abstract}

\maketitle

\section{Introduction}
\label{S:Introduction}

\subsection{Overview}
\label{SS:Overview}
One of the major breakthroughs in the theory of quantum groups in the last decade was the main quantization result from \cite{ESS}, the general explicit quantization of all classical dynamical $r$-matrices which fit Schiffmann's classification \cite{Schi}. This complemented the categorical quantization results of Etingof-Kazhdan \cite{EK} and provided a fully constructive method to quantize a given $r$-matrix. 

In this note we initiate an analogous program of constructing explicit quantizations 
in the context of Lie superalgebras. In particular we explicitly quantize zero-weight super dynamical $r$-matrices with zero coupling constant. We next discuss the classification problem for super dynamical $R$-matrices and provide some partial answers.  Then we use our results to weigh in on the question of what the correct graded analogue should be for the Hecke condition. Thus the note overall contributes to the theory of super quantum groups, which is still widely incomplete. 

\subsection{Results}
\label{SS:Results}
Our main quantization result is the following theorem, proved in Section \ref{S:QuantRes}:

\begin{theorem}
\label{T:QuantRat}
Let $\h$ be a finite dimensional commutative Lie superalgebra over $\C$ and let $V$ be a finite dimensional semisimple $\h$-module whose weights make up a basis for $\h^*$. Then every super dynamical $r$-matrix ${r : \h^* \rightarrow \End(V \otimes V)}$ with zero weight and zero coupling constant, holomorphic on an open polydisc $U \subset \h^*$, can be quantized to a super dynamical $R$-matrix $R$ on $U$.  
\end{theorem}

\noindent
(See \S\S\ref{SS:Notation}-\ref{SS:SDYBENoParam} and \S\S\ref{SS:QDYBEBasics} for the relevant definitions). 
Zero-weight super dynamical $r$-matrices with no spectral parameters were classified by the author in \cite{Kar3} in a manner which generalized the analogous non-graded results of \cite{EV}. 
In the proof of the above theorem, we make extensive use of this result, as well as results from \cite{EV2}.

The quantum theme of this note is developed mostly in Section \ref{S:QuantPict}. There we briefly study the classification problem for super dynamical $R$-matrices. and then focus on the super Hecke condition. The (non-graded) Hecke condition, introduced in \cite{EV2} as a desirable property of dynamical $R$-matrices, is a quantum analogue of the generalized unitarity condition. We proposed a super version of it in our \cite{Kar4}. In Section \ref{S:Conclusion}, we use our work here and some other considerations to weigh in on this issue of the correct super version. 

The extension to the graded world of the general constructive quantization \cite{ESS} of all classical dynamical $r$-matrices which fit Schiffmann's classification \cite{Schi} is still an open problem, and work on it is still ongoing \cite{GK}. Part of the difficulty comes from the fact that there is not yet a complete classification result analogous to \cite{Schi}; see \cite{Kar1, Kar2, Kar3, Kar4} for partial results and counterexamples in this direction.

\subsection{The organization of this note}
\label{SS:Organization}
This note is organized as follows: In Section \ref{S:Definitions} we provide the basic definitions and summarize the result from \cite{Kar3} that we will need. 
In Section \ref{S:QuantRes}, we prove Theorem \ref{T:QuantRat}. 
In Section \ref{S:QuantPict}, we sketch the development of a super analogue for the classification of super dynamical $R$-matrices given in \cite{EV2}. Section \ref{S:Conclusion} concludes the note with a discussion of the implications of our work to the problem of determining the correct way to superize the Hecke condition.

\section{Definitions and relevant earlier results}
\label{S:Definitions}

\subsection{Basic notation and terminology}
\label{SS:Notation}
Let $\g$ be a simple Lie superalgebra with non-degenerate Killing
form ${(\cdot \;,\cdot)}$. Let ${\h \subset \g}$ be a Cartan
subsuperalgebra, and let $\Delta \subset \h^*$ be the set of
roots associated to $\h$. Fix a set of simple roots $\Gamma$ or
equivalently a Borel $\b$. We will say that a set $X \subset \Delta$ of roots of $\g$ is {\em closed} if it satisfies the following:
\begin{enumerate}
\item If $\alpha, \beta \in X$ and $\alpha + \beta$ is a root,
then $\alpha + \beta \in X$, and
\item If $\alpha \in X$, then $-\alpha \in X.$
\end{enumerate}

For any positive root $\alpha$ fix $e_{\alpha} \in \g_{\alpha}$
and pick ${e_{-\alpha} \in \g_{-\alpha}}$ dual to $e_{\alpha}$
i.e. 
\[ (e_{\alpha}, e_{-\alpha}) = 1 \textmd{ for all } \alpha \in
\Delta^+.\]
\noindent
Note that we can do this uniquely up to scalars because all the $\g_{\alpha}$ are one-dimensional,
(which follows from the nondegeneracy of the Killing form).
Define:
\begin{equation} 
\label{defofA}
A_{\alpha} = \left\{ \begin{array}{cl} 
(-1)^{|\alpha|} & \textmd{if } \alpha  \textmd{ is positive} \\
1 & \textmd{if } \alpha \textmd{ is negative}
\end{array} \right. 
\end{equation}
\noindent
It is easy to see that $A_{-\alpha} = (-1)^{|\alpha|}A_{\alpha}$.
We can use $A_{\alpha}$ for instance to write the duals of our
basis vectors in terms of one another:
\[ e_{\alpha}^* = A_{-\alpha} e_{-\alpha}\]
\noindent
or equivalently:
\[ (e_{\alpha},e_{-\alpha}) = A_{-\alpha}. \]
Finally let $\Omega$ be the quadratic Casimir element, i.e.\ the element of ${\g \otimes \g}$ corresponding to the Killing form. 

The {\em super twist map} ${T_s : V \otimes V  \rightarrow V
\otimes V}$ is defined on the homogeneous elements of a given super
vector space $V$ as  
$${T_s (a\otimes b) = (-1)^{|a||b|}b\otimes a}.$$
Similarly the {\em super symmetrizing map} ${\alt : V \otimes V \otimes V  \rightarrow V \otimes V \otimes V}$ is defined on homogeneous elements by: 
\[ \alt(a\otimes b \otimes c) = a \otimes b \otimes c  +
(-1)^{|a|(|b|+|c|)} b \otimes c \otimes a + 
(-1)^{|c|(|a|+|b|)} c \otimes a \otimes b. \]

\subsection{The classical dynamical Yang-Baxter equation}
\label{SS:SDYBENoParam}
The \emph{classical dynamical Yang-Baxter equation} for a meromorphic function ${r:\h^*
\rightarrow \g \otimes \g}$ is the equation: 
\begin{equation}
\label{E:CDYBE}
\alt(dr) + 
[r^{12},r^{13}] + [r^{12},r^{23}] + [r^{13},r^{23}] = 0.
\end{equation}
Here, for a fixed (even) basis $\{x_i\}$ for $\h$, the differential of $r$ is defined as:
\[ \begin{matrix}
dr &:& \h^* &\longrightarrow& \g \otimes \g \otimes \g  \\
&& \lambda &\longmapsto& \sum_i x_i \otimes \frac{\partial
r}{\partial x_i} (\lambda)  
\end{matrix} \]
\noindent
Thus, we can see that for $r = \sum
{r}_{(1)} \otimes {r}_{(2)}$, $\alt(dr)$ may be rewritten as:
\[ \sum_i x_i^{(1)} \left(\frac{\partial r }{\partial x_i}
\right)^{(23)} +
\sum_i x_i^{(2)} \left(\frac{\partial r }{\partial x_i}
\right)^{(31)} 
+ \sum_i (-1)^{|{r}_{(1)}| |{r}_{(2)}|} x_i^{(3)} \left(
\frac{\partial r}{\partial x_i}\right)^{(12)}.\]

We will say that a meromorphic function ${r : \h^* \rightarrow \g
\otimes \g}$ is a \emph{super dynamical $r$-matrix with coupling
constant} $\epsilon$ if it is a solution to Equation
\eqref{E:CDYBE} and satisfies the {\em generalized unitarity
condition}:
\begin{equation}
\label{dynamicalunitarity} 
r(\lambda) + T_s(r)(\lambda) = \epsilon \Omega.
\end{equation}
A super dynamical $r$-matrix $r$ satisfies the \emph{zero weight
condition} if:
\[ [h\otimes 1 + 1 \otimes h, r(\lambda)] = 0 \textmd{ for all }
h \in \h, \lambda \in \h^*.\]

\subsection{Classification of super dynamical $r$-matrices of zero weight}
\label{SS:Classify}
In \cite{Kar3} we proved:
\begin{theorem}
\label{0couple0weighttheorem}
Let $\g$ be a simple Lie superalgebra with non-degenerate
Killing form ${(\cdot \; , \cdot )}$, ${\h \subset \g}$ a Cartan
subsuperalgebra, and $\Delta \subset \h^*$ the set of
roots associated to $\h$.
\begin{enumerate}
\item Let $X$ be a closed subset of the set of roots $\Delta$ of $\g$. Let $\nu \in \h^*$, and let $D = \sum_{i<j} D_{ij} dx_i \wedge
dx_j$ be a closed meromorphic $2$-form on $\h^*$. If we set
$D_{ij} = -D_{ji}$ for $i \ge j$, then the meromorphic
function:   
\begin{equation}
\label{rfor00thm}
r(\lambda) =
\sum_{i,j=1}^N D_{ij}(\lambda) x_i \otimes x_j  +  \sum_{\alpha
\in X} \frac{A_{\alpha}}{(\alpha,
\lambda-\nu)} e_{\alpha} \otimes e_{-\alpha} 
\end{equation}
\noindent
is a super dynamical $r$-matrix with zero weight and zero
coupling constant.
\item Any super dynamical $r$-matrix with zero weight and zero
coupling constant is of this form. 
\end{enumerate}
\end{theorem}

We further proved that 
there are exactly two types of zero-weight solutions to Equation \eqref{E:CDYBE} satisfying the generalized unitarity condition: the rational ones (solutions of the form given by Equation \eqref{rfor00thm}), with zero coupling constant, and the trigonometric ones, 
with a nonzero coupling constant. In fact we explicitly described the general form of the latter, but we will not need that result here.

\section{Quantization of zero weight $r$-matrices}
\label{S:QuantRes}

In this section we prove Theorem \ref{T:QuantRat}.  

\subsection{The quantum dynamical Yang-Baxter equation}
\label{SS:QDYBEBasics}
Let $\h$ be a finite dimensional commutative Lie superalgebra over $\C$, $V$ a finite dimensional super vector space over $\C$ with a diagonal(izable) $\h$ action, and let $V = \oplus_{\omega \in \h^*} V[{\omega}]$ be $V$'s $\h$-weight decomposition. In other words, for every $v \in V[{\omega}]$ and $x \in \h$, we have $x \cdot v = \omega(x)v$.

In this context, the {\em quantum dynamical Yang-Baxter equation with step $\gamma$} for a function ${R : \h^* \rightarrow \End(V \otimes V)}$ is the equation:
\begin{equation}
\label{E:QDYBE}
R^{12}(\lambda-\gamma h^{(3)})R^{13}(\lambda)R^{23}(\lambda-\gamma h^{(1)})= R^{23}(\lambda)R^{13}(\lambda-\gamma h^{(2)})R^{12}(\lambda).
\end{equation}
Here the operator $R^{ij}$ is interpreted to be acting nontrivially on the $i$th and the $j$th components of a given $3$-tensor, and the notation $h^{(k)}$ is to be replaced by the weight of the $k$th component of the same. For instance ${R^{12)}(\lambda-\gamma h^{(3)})(v_1 \otimes v_2 \otimes v_3)} = {\left (R(\lambda-\gamma \omega_3)(v_1 \otimes v_2) \right ) \otimes v_3}$ whenever $v_3 \in V[{\omega_3}]$.

We will say that an invertible function ${R : \h^* \rightarrow \End(V \otimes V)}$ is a {\em super dynamical $R$-matrix} if it is a solution to Equation \eqref{E:QDYBE} and satisfies the {\em zero weight condition}:
\[ [h\otimes 1 + 1 \otimes h, R(\lambda)] = 0 \textmd{ for all }
h \in \h, \lambda \in \h^*.\]

\subsection{The quantization problem}
\label{SS:QuantProb}
Let ${R_{\gamma} : \h^* \rightarrow \End(V \otimes V)}$ be a smooth family of solutions to Equation \eqref{E:QDYBE} such that:
\[ R_{\gamma}(\lambda) = 1 -\gamma r(\lambda) + O(\gamma^2).\]
Then the function $r(\lambda)$ satisfies Equation \eqref{E:CDYBE} and is called the {\em semi-classical limit} of $R_{\gamma}(\lambda)$. In the same setup $R_{\gamma}(\lambda)$ is called a {\em quantization} of $r(\lambda)$.

Alternatively we can begin with a super dynamical $r$-matrix $r : \h^* \rightarrow \End(V \otimes V)$ defined on an open subset $U$ of $\h^*$. We then call $r$ {\em quantizable} if there is a power series in $\gamma$ of the form:
\[ R_{\gamma}(\lambda) = 1 -\gamma r(\lambda) + \sum_{n=2}^{\infty} \gamma^nr_n(\lambda)\]
satisfying Equation \eqref{E:QDYBE}.
The {\em quantization problem} for us must now be obvious: Given a super dynamical $r$-matrix construct a power series $R_{\gamma}(\lambda)$ of the form above (or prove the impossibility of such a construction).

\subsection{Multiplicative forms}
\label{SS:Forms}

In the following we will make use of 
multiplicative $k$-forms 
a la \cite[\S\S 1.4]{EV2}. We now briefly recall some of the relevant constructions to keep our paper self-contained.

Let $V = V_{\overline{0}} \oplus V_{\overline{1}}$ be a super vector space with a homogeneous linear coordinate system $\lambda_1, \cdots, \lambda_N$. We define a {\em multiplicative $k$-form} on $V$ to be a collection:
\[ \phi = \{\phi_{i_1,\dots,i_k}(\lambda_1, \cdots, \lambda_N) \}\]
of meromorphic functions, where the ordered $k$-tuples $(i_1,\dots,i_k)$ run through all $k$-element subsets of $\{1,\cdots, N\}$, and we require that:
\[ \phi_{\tau(I)}\phi_I = 1 \]
whenever $I = (i_1,\dots,i_k)$ is some ordered $k$-tuple and $\tau(I)$ is a transposition $(i_si_{s+1})$ switching the consecutive indices $i_s,i_{s+1}$ for some $1 \le s < k$. Let $\Omega^k(V) = \Omega^k$ be the set of all multiplicative $k$-forms on $V$. There is a natural abelian group structure on $\Omega^k$.

Now we fix a complex number $\gamma$. We define, for each $i = 1, \cdots, N$, an operator $\delta_i$ on the space of all meromorphic functions on the $N$ variables $\lambda_1, \cdots, \lambda_N$:
\begin{equation*}
\delta_i : f(\lambda_1, \cdots, \lambda_N) \longmapsto 
f(\lambda_1, \cdots, \lambda_N) / f(\lambda_1, \cdots, \lambda_i-\gamma, \cdots,  \lambda_N).
\end{equation*}
We next define an operator $d_{\gamma} : \Omega^k \rightarrow \Omega^{k+1}$ mapping $\phi$ to $d_{\gamma}\phi$ given by:
\begin{equation*}
(d_{\gamma}\phi)_{i_1,\dots,i_{k+1}}(\lambda_1, \cdots, \lambda_N) = \prod_{s=1}^{k+1} \left ( \delta_{i_s} \phi_{i_1,\dots,i_{s-1}, i_{s+1}, \dots, i_{k+1}}(\lambda_1, \cdots, \lambda_N)  \right )^{(-1)^{s+1}}.
\end{equation*}
A multiplicative $k$-form $\phi$ is {\it $\gamma$-closed} if $d_{\gamma}\phi = 0$. Obviously, $d^2_{\gamma} = 0$ because the zero element of $\Omega^k$ is the form $\{\phi_{i_1,\dots,i_k}(\lambda_1, \cdots, \lambda_N)\equiv 1 \}$.

Still following \cite[\S\S 1.4]{EV2}, we say that a smooth family $\phi(\gamma) =
\{\phi_{i_1,i_2,\dots,i_k}(\lambda_1, \lambda_2, \cdots, \lambda_N,\gamma) \}$ of multiplicative $k$-forms with 
\[ \phi_I(\lambda,\gamma) = 1 - \gamma C_I(\lambda)+O(\gamma^2)  \text{ for each } I =  (i_1,i_2,\dots,i_k)\]
is a {\em quantization} of the differential form 
\[C = \sum_{i_1 < i_2 < \cdots < i_k} C_{i_1,i_2,\dots,i_k}(\lambda) \, dx_{i_1}\wedge dx_{i_2} \wedge \cdots \wedge dx_{i_k}.\]
Conversely we will say that a differential form $C$ given as above is {\em quantizable} if there exists a power series in $\gamma$:
\[ \phi_I(\lambda,\gamma) = 1 - \gamma C_I(\lambda)+\sum_{n=2}^{\infty} \gamma^n C_{n;I}(\lambda) \text{ for each } I =  (i_1,i_2,\dots,i_k)\]
convergent for small $|\gamma|$ and fixed $\lambda \in U$, where $U$ is an open polydisc in $\C^N$, in such a way that $\{\phi_{i_1,i_2,\dots,i_k}(\lambda_1, \lambda_2, \cdots, \lambda_N,\gamma) \}$ is a multiplicative $k$-form. 

Here is Lemma 1.1 from \cite{EV2}:
\begin{lemma}
\label{L:Lemma1.1}
Every closed holomorphic differential $k$-form $C$ defined on an open polydisc is quantizable to a holomorphic multiplicative closed $k$-form $\phi(\gamma)$.
\end{lemma}
\noindent
The proof is included in \cite{EV2} and will not be repeated here.

\subsection{$R$-matrices of $\gl(m,n)$ type}
\label{SS:glmntype}
Now let $\h$ be a finite dimensional commutative Lie superalgebra over $\C$ and let $V = V_{\overline{0}} \oplus V_{\overline{1}}$ be a finite dimensional semisimple $\h$-module whose weights $W= \{\omega_1, \omega_2, \cdots, \omega_N\}$ make up a basis for $\h^*$. We label the elements of the dual basis for $\h$ by $x_i$; clearly the $x_i$ are all even, and $\dim_{\C} V = \dim_{\C} \h = N$. Let $\{v_1, v_2, \cdots, v_N\}$ be an $\h$-eigenbasis for $V$ with $x_i v_j = \delta_{ij}v_j$. By relabeling as needed, we can assume that $v_1, v_2, \cdots v_m$ is a basis for $V_{\overline{0}}$, the even part of $V$, while $v_{m+1},v_{m+2} \cdots, v_{N}$ is a basis for $V_{\overline{1}}$, the odd part of $V$. Let $n = N-m$. We will say that a super dynamical $R$-matrix $R : \h^* \rightarrow \End(V \otimes V)$ for such $\h$ and $V$ is an {\em $R$-matrix of $\gl(m,n)$ type}. The super dynamical $R$-matrices in this paper will all be of this kind unless explicitly noted otherwise. 

In this setup $V \otimes V$ has the following weight decomposition:
\begin{equation}
\label{E:WeightDecomp}
V \otimes V = \left ( \bigoplus_{i=1}^N V_{ii} \right ) \oplus \left ( \bigoplus_{i <j} V_{ij}  \right ). 
\end{equation}
Here $V_{ii} = \C (v_i \otimes v_i)$ and $V_{ij} = \C (v_i \otimes v_j) \oplus \C (v_j \otimes v_i)$. It is clear that $V_{ii}$ will always belong to the even part of $V \otimes V$ while $V_{ij}$ may be even or odd. In particular, if exactly one of $v_i$ and $v_j$ is odd, then $V_{ij}$ will be odd, otherwise it will be even.  In other words, if we introduce the notation 
\[ \sigma(i) = \begin{cases}
0 & \text{ if } i \le m \\
1 & \text{ if } i > m,
\end{cases}\]
then $V_{ij}$ is odd if and only if $\sigma(i)+\sigma(j) = 1$. 

We can introduce a basis $\{E_{ij} \, \vert 1 \le i,j\le N\}$  for $\End(V)$ by setting $E_{ij} v_k = \delta_{jk} v_i$. Recall that we require our dynamical $R$-matrices to satisfy the zero weight condition. Then we can write any $R$-matrix $R : \h^* \rightarrow \End(V\otimes V)$ of $\gl(m,n)$ type in the form:
\[ R(\lambda) = \sum_{i,j=1}^N \alpha_{ij}(\lambda) E_{ii} \otimes E_{jj} + \sum_{i \neq j} \beta_{ij}(\lambda) E_{ji}\otimes E_{ij} \] 
for some meromorphic functions $\alpha_{ij}, \beta_{ij} : \h^* \rightarrow \C$. 

\subsection{Gauge transformations for super dynamical $r$-matrices}
\label{SS:Gauge}

Before we can prove Theorem \ref{T:QuantRat},
we will need to simplify the expression \eqref{rfor00thm}. In order to do that we first discuss briefly the appropriate gauge transformations for super dynamical $r$-matrices of the form $r : \h^* \rightarrow \End(V \otimes V)$. Note that if we assume the setup of \S\S\ref{SS:glmntype}, then we can use the $\{E_{ij}\}$ basis for $\End(V)$. Then the zero weight condition on $r$ implies that $r$ has to be in the form:
\[ r(\lambda) = \sum_{i,j=1}^N \alpha_{ij}(\lambda) E_{ii} \otimes E_{jj} + \sum_{i \neq j} \beta_{ij}(\lambda) E_{ji}\otimes E_{ij} \] 
for some meromorphic functions $\alpha_{ij}, \beta_{ij} : \h^* \rightarrow \C$. 

The following is a list of the gauge transformations for such $r$ which we will need in the rest of this note (cf.\ \cite{EV, EV2}):
\begin{enumerate}\addtolength{\itemsep}{0.4\baselineskip}
\item The transformation:
\[ r(\lambda) \longmapsto r(\lambda) + \sum_{i,j=1}^N D_{ij}(\lambda) E_{ii} \otimes E_{jj} \]
for some closed meromorphic differential $2$-form $D = \sum_{i<j} D_{ij} dx_i \wedge dx_j$ on $\h^*$. $D_{ij}$ is then extended to all $i,j$ by setting
$D_{ij} = -D_{ji}$ for $i \ge j$. ($D_{ii} = 0$ for each $i$.)
\item The transformation:
\[ r(\lambda) \longmapsto r(\lambda+\mu) \]
for $\mu \in \h^*$.
\item The transformation:
\[ r(\lambda) \longmapsto c r(c\lambda) \]
for a nonzero complex number $c \in \C$.
\item The transformation:
\[ r(\lambda) \longmapsto (\tau \otimes \tau) r(\tau^{-1} \cdot \lambda) (\tau^{-1} \otimes \tau^{-1}) \]
for some permutation $\tau \in S_N$
of the coordinates in $\h^*$ and $V$.
\item The transformation: 
\[ r(\lambda) \longmapsto r(\lambda) + c\operatorname{Id}\]
for a nonzero complex number $c \in \C$.
\end{enumerate}
Each of these transformations corresponds to a specific quantum gauge transformation allowed for super dynamical $R$-matrices (cf.\ \cite{EV2}). 
We will briefly study  
these quantum gauge transformations 
in \S\S\ref{SS:SuperQuantGauge}, in the context of the classification problem for super dynamical $R$-matrices.  

It is easy to show that the transformations (1-5) map a given super dynamical $r$-matrix to another. We omit the proofs here since they are straightforward modifications of those in \cite{EV}. We will say that two super dynamical $r$-matrices are {\em gauge equivalent} (or simply {\em equivalent} when the context is unambiguous) if one can be obtained from the other by a sequence of gauge transformations.

We can now simplify the expression in Theorem \ref{0couple0weighttheorem} using the above. Let $X \subset \{ 1, 2, \cdots, N\}$ be a subset of indices and write it as a disjoint union of subintervals $X = X_1 \sqcup X_2 \sqcup \cdots \sqcup X_n$. In other words, every subinterval $X_k$ should be of the form $X_k = [i_k,i_{k+1}, i_{k+2}, \cdots, j_k]$, and $j_k < i_{k+1}$ for each $k$. Define:
\begin{equation*} 
A_{ij} =  \left\{ \begin{array}{cl} 
(-1)^{\sigma(i)+\sigma(j)} &\text{ if } i < j, \\
1 & \text{ if } i > j,
\end{array} \right. 
\end{equation*}
(cf.\ Equation \eqref{defofA}). Now applying the above transformations 
and using the $\{E_{ij}\}$ basis, we can show that the super dynamical $r$-matrix in Equation \eqref{rfor00thm} is (gauge-)equivalent to:
\[ r_{\text{rat}}(\lambda) =  \sum_{k=1}^n \left ( \sum_{i,j \in X_k, i\neq j} \frac{A_{ij}}{\lambda_{ij}} E_{ij} \otimes (E_{ij})^* \right ). \]
Since $(E_{ij})^* = A_{ji}(-1)^{\sigma(i)}E_{ji}$, this further reduces to:
\begin{equation}
\label{E:rRatForm}
r_{\text{rat}}(\lambda) =  \sum_{k=1}^n \left ( \sum_{i,j \in X_k, i\neq j} \frac{(-1)^{\sigma(j)}}{\lambda_{ij}} E_{ij} \otimes E_{ji} \right ). 
\end{equation}

\subsection{The construction}
\label{SS:Constructions}
We are finally ready to construct the quantization necessary for Theorem \ref{T:QuantRat}. Let $\h$ and $V$ be as in \S\S\ref{SS:glmntype}. We will once again use the basis $\{E_{ij}\}$ for $\End(V)$ and we will write $X \subset \{ 1, 2, \cdots, N\}$ as a disjoint union of subintervals $X = X_1 \sqcup X_2 \sqcup \cdots \sqcup X_n$

Consider:
\[ R_{\text{rat}}(\lambda,\gamma) = \text{Id} + \sum_{k=1}^n \, \sum_{i,j \in X_k,i\neq j} \frac{\gamma}{\lambda_{ij}} \left ( E_{ii} \otimes E_{jj} + (-1)^{\sigma(i)} E_{ji}\otimes E_{ij} \right )\]
Then $R_{\text{rat}}(\lambda,\gamma)$ satisfies Equation \eqref{E:QDYBE} (cf.\ Theorem \ref{T:ClassifyR1}), and its semi-classical limit is:
\[r^{\prime}_{\text{rat}}(\lambda) = \sum_{k=1}^n \, \sum_{i,j \in X_k,i\neq j} \frac{-1}{\lambda_{ij}} \left ( E_{ii} \otimes E_{jj} + (-1)^{\sigma(i)} E_{ji}\otimes E_{ij} \right ).\]
Using the gauge transformation of type (1) with the closed form:
\[ D =\sum_{k=1}^n \, \sum_{i,j \in X_k,i< j} D_{ij} dx_i \wedge dx_j = \sum_{k=1}^n \, \sum_{i,j \in X_k,i< j} \frac{-1}{\lambda_{ij}}  dx_{i} \wedge dx_{j},  \]
we can show that $r^{\prime}_{\text{rat}}$ is (gauge-)equivalent to the $r_{\text{rat}}$ of Equation \eqref{E:rRatForm}.
Together with Lemma \ref{L:Lemma1.1} this proves Theorem \ref{T:QuantRat}. \qed

\section{The Quantum Picture}
\label{S:QuantPict}

In this section we define {\it the super Hecke condition} (\S\S\ref{SS:SuperHecke}) which is a generalized unitarity condition. Using this notion, we state and prove (\S\S\ref{SS:TheoremsStated}) a theorem in the spirit of Theorem 1.2 of \cite{EV2}. This is a result that provides a partial classification of all super dynamical $R$-matrices satisfying the super Hecke condition. It turns out that the super Hecke condition encodes the constraint on the coupling constant in the classical case.

\subsection{Some initial computations}
\label{SS:QuantumComputations}

Let $\h$ and $V$ be as in \S\S\ref{SS:glmntype}. We will once again use the basis $\{E_{ij}\}$ for $\End(V)$ and throughout this section we will once again restrict ourselves to the study of $R$-matrices of $\gl(m,n)$ type. Recall that this means, in particular, that the super vector space $V \otimes V$ has the weight decomposition given in \eqref{E:WeightDecomp}. 

More specifically, a super dynamical $R$-matrix $R : \h^* \rightarrow \End(V\otimes V)$ of $\gl(m,n)$ type can be written in the form:
\[ R(\lambda) = \sum_{i,j=1}^N \alpha_{ij}(\lambda) E_{ii} \otimes E_{jj} + \sum_{i \neq j} \beta_{ij}(\lambda) E_{ji}\otimes E_{ij} \] 
for some meromorphic functions $\alpha_{ij}, \beta_{ij} : \h^* \rightarrow \C$. 
If we now assume for simplicity (and for other reasons which will become clearer in \S\S\ref{SS:SuperHecke}) that our super dynamical $R$-matrices all satisfy $\alpha_{ii} = 1$ for all $i$, we can rewrite the above as:
\begin{equation}
\label{E:FormofR}
R(\lambda) = \sum_{i=1}^N E_{ii} \otimes E_{ii} + \sum_{i \neq j} \alpha_{ij}(\lambda) E_{ii} \otimes E_{jj} + \sum_{i \neq j} \beta_{ij}(\lambda) E_{ji}\otimes E_{ij} \end{equation}
for some meromorphic functions $\alpha_{ij}, \beta_{ij} : \h^* \rightarrow \C$. 

In this subsection we list a few conditions on these $\alpha$ and $\beta$ functions. We limit ourselves to simply summarizing the results of necessary computations; the explicit derivations can be found in Appendix \ref{A:AlphaBeta}. 

By applying the two sides of Equation \eqref{E:QDYBE} for an $R$ of the form \eqref{E:FormofR} to a basis element $v_i \otimes v_i \otimes v_k$ of $V^{\otimes 3}$ with $i \neq k$ and setting the coefficients of like terms equal to one another, we obtain:
\begin{equation}
\label{E:coeffkii}
\alpha_{ki}(\lambda-\gamma \omega_i)\beta_{ik}(\lambda)\alpha_{ik}(\lambda-\gamma \omega_i) + (\beta_{ik}(\lambda-\gamma \omega_i))^2= \beta_{ik}(\lambda-\gamma \omega_i)
\end{equation}
and 
\begin{eqnarray}
\label{E:coeffiki}
{(-1)^{\sigma(i)+\sigma(k)}} \beta_{ki}(\lambda-\gamma \omega_i)\beta_{ik}(\lambda) \alpha_{ik}(\lambda-\gamma \omega_i) &+& \alpha_{ik}(\lambda-\gamma \omega_i)\beta_{ik}(\lambda-\gamma \omega_i) \cr
&=& \beta_{ik}(\lambda)\alpha_{ik}(\lambda-\gamma \omega_i). 
\end{eqnarray}
Note that Equation \eqref{E:coeffkii} is identical to \cite[Eqn.1.8.4]{EV2} while Equation \eqref{E:coeffiki} is a signed version of \cite[Eqn.1.8.5]{EV2}. 

Similarly we can derive the following equations by applying the two sides of Equation \eqref{E:QDYBE} to a basis element $v_i \otimes v_j \otimes v_k$ with $i,j,k$ all distinct:
\begin{equation}
\label{E:coeffijk}
\alpha_{ij}(\lambda-\gamma \omega_k)\alpha_{ik}(\lambda)\alpha_{jk}(\lambda-\gamma \omega_i) = \alpha_{jk}(\lambda)\alpha_{ik}(\lambda-\gamma \omega_j)\alpha_{ij}(\lambda)
\end{equation}
which is precisely the same as \cite[Eqn.1.8.6]{EV2};
\begin{equation}
\label{E:coeffikj}
\alpha_{ik}(\lambda-\gamma \omega_j)\alpha_{ij}(\lambda) \beta_{jk}(\lambda-\gamma \omega_i)=  \beta_{jk}(\lambda)\alpha_{ik}(\lambda-\gamma \omega_j)\alpha_{ij}(\lambda)
\end{equation}
which is precisely the same as \cite[Eqn.1.8.7]{EV2};
\begin{equation}
\label{E:coeffjik}
\beta_{ij}(\lambda-\gamma \omega_k)\alpha_{ik}(\lambda)\alpha_{jk}(\lambda-\gamma \omega_i) = \alpha_{ik}(\lambda)\alpha_{jk}(\lambda-\gamma\omega_i)\beta_{ij}(\lambda)
\end{equation}
which is precisely the same as \cite[Eqn.1.8.8]{EV2};
\begin{eqnarray}
\label{E:coeffjki}
(-1)^{\sigma(k)}\beta_{kj}(\lambda-\gamma \omega_i)\beta_{ik}(\lambda) \alpha_{jk}(\lambda-\gamma \omega_i)&+&(-1)^{\sigma(j)} \alpha_{jk}(\lambda-\gamma \omega_i)\beta_{ij}(\lambda)  \beta_{jk}(\lambda-\gamma \omega_i)\cr
&=&(-1)^{\sigma(i)} \beta_{ik}(\lambda)\alpha_{jk}(\lambda-\gamma\omega_i)\beta_{ij}(\lambda)
\end{eqnarray}
which is a signed analogue of \cite[Eqn.1.8.9]{EV2};
\begin{eqnarray}
\label{E:coeffkji}
\alpha_{kj}(\lambda-\gamma \omega_i)\beta_{ik}(\lambda)\alpha_{jk}(\lambda-\gamma \omega_i) &+& \beta_{jk}(\lambda-\gamma \omega_i)\beta_{ij}(\lambda) \beta_{jk}(\lambda-\gamma \omega_i) =\cr
\alpha_{ji}(\lambda)\beta_{ik}(\lambda-\gamma \omega_j)\alpha_{ij}(\lambda) &+&(-1)^{\sigma(i)+\sigma(j)}\beta_{ij}(\lambda)\beta_{jk}(\lambda-\gamma\omega_i)\beta_{ij}(\lambda)\end{eqnarray}
which is a signed analogue of \cite[Eqn.1.8.10]{EV2}; and
\begin{eqnarray}
\label{E:coeffkij}
\beta_{ik}(\lambda-\gamma \omega_j)\alpha_{ij}(\lambda)\beta_{jk}(\lambda-\gamma \omega_i)&=&\cr
 \beta_{ji}(\lambda)\beta_{ik}(\lambda-\gamma \omega_j)\alpha_{ij}(\lambda) &+&\alpha_{ij}(\lambda)\beta_{jk}(\lambda-\gamma\omega_i)\beta_{ij}(\lambda)\end{eqnarray}
which is precisely the same as \cite[Eqn.1.8.11]{EV2}.

\subsection{The Super Hecke Condition}
\label{SS:SuperHecke}

Let $p \neq -q$ be two complex numbers. Set $\check{R} = P_s R$ where $P_s \in \End(V \otimes V)$ is the element corresponding to $T_s$.
In a way analogous to \cite{EV2} we will say that a function $R : \h^* \rightarrow \End(V \otimes V)$ satisfies the {\em strong super Hecke condition} if it has the following properties:
\begin{enumerate}\addtolength{\itemsep}{0.4\baselineskip}
\item The function preserves the weight decomposition given in \eqref{E:WeightDecomp}.
\item For any $i =1, 2, \cdots, N$, and $\lambda \in \h^*$, $\check{R}(\lambda)(v_i \otimes v_i) = p (v_i \otimes v_i)$.
\item For any $i \neq j$, and $\lambda \in \h^*$,  the operator $\check{R}(\lambda)$ restricted to $V_{ij}$ has eigenvalues $(-1)^{\sigma(i)+\sigma(j)}p$ and $-(-1)^{\sigma(i)+\sigma(j)}q$. 
\end{enumerate}

A function $R : \h^* \rightarrow \End(V \otimes V)$ satisfies the {\em weak super Hecke condition} if it has the following properties  (cf\ \cite[Eq.1.3.6]{EV2}) :
\begin{enumerate}\addtolength{\itemsep}{0.4\baselineskip}
\item The function preserves the weight decomposition given in \eqref{E:WeightDecomp}.
\item For any $\lambda \in \h^*$ and $i,j \le N$, $(\check{R}(\lambda) -(-1)^{\sigma(i)+\sigma(j)} p)(\check{R}(\lambda)+(-1)^{\sigma(i)+\sigma(j)}q)= 0$ when restricted to $V_{ij}$.
\end{enumerate}

Just as in the non-graded case these two properties are intimately related. In fact whenever a continuous family $R_t : \h^* \rightarrow \End(V\otimes V)$, $t \in [0,1]$, of meromorphic functions, analytic for $0 < t < 1$ and $R_0 = \operatorname{Id}$, satisfies the weak super Hecke condition for all $t$,  then $R_t$ satisfies the strong super Hecke condition as well. 
Hence we will simply assume that $R$ satisfies both whenever we say that $R$ satisfies the super Hecke condition.

Now we consider a super dynamical $R$-matrix $R(\lambda)$ with step $\gamma=1$ which satisfies the super Hecke condition with $p=1$ and $q$ arbitrary. Then we can see that $\alpha_{ii} = 1$ and $R$ has the form given by Equation \eqref{E:FormofR}. Furthermore, whenever $i \neq j$, we have:
\begin{equation}
\label{E:trace}
(-1)^{\sigma(i)}\beta_{ij}(\lambda) + (-1)^{\sigma(j)}\beta_{ji}(\lambda) = (-1)^{\sigma(i)+\sigma(j)}(1-q)
\end{equation}
and 
\begin{equation}
\label{E:determinant}
(-1)^{\sigma(i)+\sigma(j)}\beta_{ij}(\lambda)\beta_{ji}(\lambda)-\alpha_{ij}(\lambda)\alpha_{ji}(\lambda)=-q
\end{equation}
obtained from the trace and determinant of $\check{R}$ on $V_{ij}$. Note that these are signed versions of \cite[Eqn.1.8.2]{EV2}  and \cite[Eqn.1.8.3]{EV2}.

At this point it is easy to notice that if $i\neq j$, then assuming $\alpha_{ij} \equiv 0$ implies that 
\[ \beta_{ij}(\lambda)\beta_{ji}(\lambda) = -(-1)^{\sigma(i)+\sigma(j)}q \]
by Equation \eqref{E:determinant}, and Equation \eqref{E:coeffkii} gives us:
\[ (\beta_{ij}(\lambda))^2= \beta_{ij}(\lambda) \text{ and } (\beta_{ji}(\lambda))^2= \beta_{ji}(\lambda). \]
These then contradict with Equation \eqref{E:trace}. Therefore $\alpha_{ij}$ cannot be identically zero. Similarly we can show that 
\begin{equation}
\label{E:alphasbetas}
\alpha_{ij}(\lambda)\alpha_{ji}(\lambda) = ((-1)^{\sigma(i)}\beta_{ij}(\lambda)+(-1)^{\sigma(i)+\sigma(j)}q) ((-1)^{\sigma(j)}\beta_{ji}(\lambda)+(-1)^{\sigma(i)+\sigma(j)}q)
\end{equation}
and therefore the quantity $(-1)^{\sigma(i)}\beta_{ij}(\lambda)+(-1)^{\sigma(i)+\sigma(j)}q$ is also not identically zero.

Finally we consider a super dynamical $R$-matrix $R(\lambda)$ of the form \eqref{E:FormofR} with step $\gamma=1$, and assume that $R(\lambda)$ satisfies the super Hecke property with Hecke parameters $p=1$ and $q$. Then the collection of functions:
\begin{equation}
\label{E:phiDef}
\phi = \{ \phi_{ij}(\lambda) \} \text{ where } \phi_{ij}(\lambda) = \frac{(-1)^{\sigma(i)}\beta_{ij}(\lambda)+(-1)^{\sigma(i)+\sigma(j)}q}{\alpha_{ij}(\lambda)} \text{ for } i\neq j 
\end{equation}
is a $\gamma$-closed multiplicative $2$-form with $\gamma = 1$. This follows from our earlier computations and in particular from Equation \eqref{E:alphasbetas}; just as in \cite{EV2}, Equations \eqref{E:coeffijk} and \eqref{E:coeffikj} are used to show that $d_{\gamma}\phi = 0$. We will use this $\phi$ in the next subsection.

\subsection{Gauge transformations for super dynamical $R$-matrices}
\label{SS:SuperQuantGauge}

Let us now assume that we have a super dynamical $R$-matrix of $\gl(m,n)$ type and we write it in the form given by Equation \eqref{E:FormofR}. The following is a list of the gauge transformations for such $R(\lambda)$ which we will need in the rest of this note (cf.\ \cite[\S\S1.4]{EV2}):
\begin{enumerate}
\item The transformation:
\[ R(\lambda) \longmapsto  \sum_{i=1}^N E_{ii} \otimes E_{ii} + \sum_{i \neq j} \phi_{ij}(\lambda) \alpha_{ij}(\lambda) E_{ii} \otimes E_{jj} + \sum_{i \neq j} \beta_{ij}(\lambda) E_{ji}\otimes E_{ij} \]
for some meromorphic s-multiplicative $\gamma$-closed multiplicative $2$-form $\{\phi_{ij}\}$ on $\h^*$.
\item The transformation:
\[  R(\lambda) \longmapsto (\tau \otimes \tau) R(\tau^{-1} \cdot \lambda) (\tau^{-1} \otimes \tau^{-1}) \] 
for some permutation $\tau \in S_N$
of the coordinates in $\h^*$ and $V$.
\item The transformation:
\[ R(\lambda) \longmapsto cR(\lambda) \]
for a nonzero complex number $c \in \C$.
\item The transformation:
\[ R(\lambda) \longmapsto R(c\lambda+\mu) \]
for a nonzero complex number $c \in \C$ and an element $\mu \in \h^*$.
\end{enumerate}

It is easy to see that transformations of type (1-3) transform a super dynamical $R$-matrix with step $\gamma$ to another one with step $\gamma$. In particular it suffices to check that the relevant equations in \S\S\ref{SS:QuantumComputations} and \S\S\ref{SS:SuperHecke} for $\alpha_{ij}(\lambda)$ and $\beta_{ij}(\lambda)$ are invariant with respect to them. 
 Transformations of type (4) modify the step $\gamma$ to $\gamma/c$. 
 
In all cases the super Hecke property is preserved. While the transformations of type (3) modify the relevant Hecke parameters, the rest preserve them. Moreover, any super dynamical $R$-matrix $R(\lambda)$ of Hecke type can be shown to be (gauge-)equivalent to a super dynamical $R$-matrix $R(\lambda)$ with step $\gamma=1$ which satisfies the super Hecke condition with $p=1$ and $q$ arbitrary. This requires simply gauge transformations of types (3) and (4). 

At this point we can specialize \eqref{E:FormofR} even further. Once again let $R(\lambda)$ be a super dynamical $R$-matrix of the form \eqref{E:FormofR} with step $\gamma$ satisfying the super Hecke condition. As justified by the above we can assume that the step $\gamma=1$ and the Hecke parameters are $p=1$ and $q$ arbitrary. Then if we apply the gauge transformation of type (1) to this $R(\lambda)$ using the reciprocal of the multiplicative $2$-form given in \eqref{E:phiDef}, we obtain a new super dynamical $R$-matrix (satisfying the super Hecke condition with the same parameters) whose coefficients now satisfy 
\begin{equation}
\label{E:alphabeta}
 {\alpha_{ij}(\lambda)} =  (-1)^{\sigma(i)}\beta_{ij}(\lambda)+(-1)^{\sigma(i)+\sigma(j)}q \text{ for } i\neq j. 
 \end{equation}

\subsection{Statement of the main quantum theorem}
\label{SS:TheoremsStated}

We are now ready to state the main result of this section:

\begin{theorem}[Classification Theorem for Equal Parameters]
\label{T:ClassifyR1}
Let $\h$ be a finite dimensional commutative Lie superalgebra over $\C$ and let $V = V_{\overline{0}} \oplus V_{\overline{1}}$ be a finite dimensional semisimple $\h$-module whose weights make up a basis for $\h^*$. Let $N = \dim_{\C} V = \dim_{\C} \h$ . 
\begin{enumerate}
\item 
Let $X \subset \{ 1, 2, \cdots, N\}$ be a subset of indices written as a disjoint union of subintervals $X = X_1 \sqcup X_2 \sqcup \cdots \sqcup X_n$. 
Fix a $\gamma$-quasiconstant $\mu : \h^* \rightarrow \h^*$ with $\gamma=1$. Define scalar meromorphic $\gamma$-quasiconstant functions $\mu_{ij} : \h^* \rightarrow \C$ by $\mu_{ij}(\lambda) = x_1(\mu(\lambda)) - x_j(\mu(\lambda))$. 
Then the meromorphic function $R_X : \h^* \rightarrow \End(V \otimes V)$ defined by:
\begin{equation*}
R_X(\lambda) = \sum_{i,j=1}^N (-1)^{\sigma(i)+\sigma(j)}E_{ii} \otimes E_{jj} + \sum_{s=1}^n \left ( \sum_{i,j \in X_s, i \neq j}  \frac{1}{\lambda_{ij} - \mu_{ij}(\lambda)} [E_{ii}\otimes E_{jj} + (-1)^{\sigma(i)}E_{ji} \otimes E_{ij}] \right ) 
\end{equation*}
is a super dynamical $R$-matrix of $\gl(m,n)$ type satisfying the super Hecke condition with $p=q=1$ and step $\gamma=1$. 
\item Every super dynamical $R$-matrix of $\gl(m,n)$ type satisfying the super Hecke condition with $p=q$ is equivalent to a super dynamical $R$-matrix of this form. 
\end{enumerate}
\end{theorem}

\subsection{Proof of Theorem \ref{T:ClassifyR1}}
\label{SS:QProof1}

Now we let $R(\lambda)$ be a super dynamical $R$-matrix satisfying the super Hecke condition with parameters $p=q$. 
As we showed in the previous subsection, we can use appropriate gauge transformations to ensure that $\gamma = p=q=1$. Then Equation \eqref{E:alphabeta} becomes:
\[  {\alpha_{ij}(\lambda)} =  (-1)^{\sigma(i)}\beta_{ij}(\lambda)+(-1)^{\sigma(i)+\sigma(j)} \text{ for } i\neq j. \]
Next look at Equation \eqref{E:coeffiki} for indices $i,j$. Clearly $\beta_{ij}(\lambda) = \beta_{ji}(\lambda) \equiv 0$ is one solution, so we assume that this is not the case. Since we showed earlier in \S\S\ref{SS:SuperHecke} that $\alpha_{ij}$ cannot be identically zero, we obtain from the two versions (for $i,i,j$ and $j,j,i$, reading the coefficients of $i,j,i$ and $j,i,j$ respectively), the following two conditions on $\beta_{ij}$:
\begin{equation}
\label{E:betaReciprocal1}
\frac{1}{\beta_{ij}(\lambda)} - \frac{1}{\beta_{ij}(\lambda-\omega_i)} = 1
 \text{ for } i\neq j, \end{equation}
and
\begin{equation}
\label{E:betaReciprocal2}
\frac{1}{\beta_{ij}(\lambda)} - \frac{1}{\beta_{ij}(\lambda-\omega_j)} = -(-1)^{\sigma(i)+\sigma(j)} \text{ for } i\neq j. \end{equation}
where we are using $(-1)^{\sigma(i)}\beta_{ij}(\lambda) + (-1)^{\sigma(j)}\beta_{ji}(\lambda)=0$ or equivalently $A_{ij}\beta_{ij}(\lambda) + A_{ji}\beta_{ji}(\lambda)=0$
(obtained from Equation \eqref{E:trace} with $q=1$).

Rewriting these equations as:
\begin{equation*}
\beta_{ij}(\lambda-\omega_i) = \frac{\beta_{ij}(\lambda)}{1-\beta_{ij}(\lambda)}
\end{equation*}
and
\begin{equation*}
\beta_{ji}(\lambda-\omega_i) = \frac{\beta_{ji}(\lambda)}{1+ (-1)^{\sigma(i)+\sigma(j)}\beta_{ji}(\lambda)}
\end{equation*}
and using the description of $\alpha_{ij}(\lambda)$ in terms of the $\beta_{ij}(\lambda)$ given above, we see that solutions $\beta_{ij}(\lambda), \beta_{ji}(\lambda)$ to the above equations will also be solutions to Equation \eqref{E:coeffkii} (cf.\ \cite[Lemma 1.4]{EV2}.

Furthermore defining 
\begin{equation*}
\mu_{ij}(\lambda) = \lambda_{ij} - \frac{(-1)^{\sigma(i)}}{\beta_{ij}(\lambda)}
\end{equation*}
we can show that $\mu_{ij}(\lambda-\omega_i) = \mu_{ij}(\lambda-\omega_j) = \mu_{ij}(\lambda)$ for all $i \neq j$.
Thus the meromorphic functions
\begin{equation*}
\beta_{ij}(\lambda) = \frac{(-1)^{\sigma(i)}}{\lambda_{ij} - \mu_{ij}(\lambda)} \text{ and } 
\beta_{ji}(\lambda) = \frac{(-1)^{\sigma(j)}}{\lambda_{ji} - \mu_{ji}(\lambda)}  
\end{equation*} 
where $\mu_{ij}(\lambda)=-\mu_{ji}(\lambda)$ and $\mu_{ij}(\lambda)$ is a meromorphic function periodic with respect to shifts of $\lambda$ by $\omega_i$ and $\omega_j$ will be solutions to Equations \eqref{E:coeffkii} and \eqref{E:coeffiki} (cf.\ \cite[Lemma 1.4]{EV2}).

Note that Equations \eqref{E:coeffikj} and \eqref{E:coeffjik} imply that the function $\beta_{ij}(\lambda)$ is periodic with respect to shifts of $\lambda$ by $\omega_k$ for all $k$ distinct from $i$ and $j$. This periodicity then holds also for $\mu_{ij}(\lambda)$. 

Next look at Equation \eqref{E:coeffjki} on functions $\beta_{ij}(\lambda)$, $\beta_{jk}(\lambda)$, $\beta_{ik}(\lambda)$, we note that if any one of these is identically zero, then at least one more has to be identically zero. This allows us to define an equivalence relation on the indices $\{1, 2, \cdots, N\}$: First assert that all $i$ are related to themselves. Then for $i \neq j$ let $i$ be related to $j$ if $\beta_{ij}(\lambda)$ is not identically zero. The symmetry property follows directly from the trace condition. 

For the equivalence relation defined above, let $Y = Y_1 \cup Y_2 \cup \cdots \cup Y_n$ be the union of all $n$ equivalence classes $Y_i$ with more than one element. If pairwise distinct $i,j,k$ do not all belong in the same equivalence class, then at least two of $\beta_{ij}$, $\beta_{jk}$, $\beta_{ik}$ will be identically zero, thus the triple will be consistent with Equation \eqref{E:coeffjki}. If all three lie in the same equivalence class, then we get:
\begin{equation*}
(-1)^{\sigma(k)}\beta_{kj}(\lambda-\omega_i)\beta_{ik}(\lambda) + (-1)^{\sigma(j)}\beta_{ij}(\lambda)  \beta_{jk}(\lambda-\omega_i) = (-1)^{\sigma(i)} \beta_{ik}(\lambda)\beta_{ij}(\lambda),
\end{equation*}
and by periodicity of $\beta_{kj}$ and $\beta_{jk}$ with respect to $\omega_i$, we reduce this further to:
\begin{equation*}
(-1)^{\sigma(k)}\beta_{kj}(\lambda)\beta_{ik}(\lambda) + (-1)^{\sigma(j)}\beta_{ij}(\lambda)\beta_{jk}(\lambda) = (-1)^{\sigma(i)} \beta_{ik}(\lambda)\beta_{ij}(\lambda).
\end{equation*}
We can rewrite this as:
\begin{align*}
(-1)^{\sigma(k)}\left ( \frac{(-1)^{\sigma(k)}}{\lambda_{kj} - \mu_{kj}(\lambda)}\right ) \left ( \frac{(-1)^{\sigma(i)}}{\lambda_{ik} - \mu_{ik}(\lambda)}\right ) 
&+ (-1)^{\sigma(j)}\left ( \frac{(-1)^{\sigma(i)}}{\lambda_{ij} - \mu_{ij}(\lambda)} \right ) \left ( \frac{(-1)^{\sigma(j)}}{\lambda_{jk} - \mu_{jk}(\lambda)} \right ) \\
&= (-1)^{\sigma(i)} \left ( \frac{(-1)^{\sigma(i)}}{\lambda_{ik} - \mu_{ik}(\lambda)}\right ) \left ( \frac{(-1)^{\sigma(i)}}{\lambda_{ij} - \mu_{ij}(\lambda)}\right ), 
\end{align*}
which, after sign cancelations, yields $\mu_{ik}(\lambda) = \mu_{ij}(\lambda)+\mu_{jk}(\lambda)$. Therefore as in the non-graded case of \cite[\S\S1.11]{EV2}, we conclude that there exists a $1$-quasiconstant meromorphic function $\mu : \h^* \rightarrow \h^*$ such that $\mu_{ij}(\lambda) = x_i(\mu(\lambda)) - x_j(\mu(\lambda))$ for all $i,j$ with $\mu_{ij}$ not identically zero, and thus 
Equation \eqref{E:coeffkji} is also satisfied.

Let $\tau$ be a permutation of $\{1, \dots, N\}$ that transforms the set $Y$ into a set $X$ which can now be written as a disjoint union of subintervals $X = X_1 \sqcup X_2 \sqcup \cdots \sqcup X_n$. In other words, every subinterval $X_k$ should be of the form $X_k = [i_k,i_{k+1}, i_{k+2}, \cdots, j_k]$, and $j_k < i_{k+1}$ for each $k$. Finally applying a gauge transformation of type (2) for this $\tau$ to the $R$-matrix $R$ will yield an $R$-matrix of the form desired. This completes the proof of Theorem \ref{T:ClassifyR1}. \qed

\section{Conclusion}
\label{S:Conclusion}

In this note we proved a quantization theorem for super dynamical $r$-matrices. More specifically we explicitly constructed quantizations for zero weight super dynamical $r$-matrices with zero coupling constant. We expect that the definitions and constructions here will also be helpful in the proof of an analogous quantization result for nonzero coupling constants, we plan to follow up on this thread in future work. 

It must be clear that quantization in this note meant finding a solution to the quantum dynamical Yang-Baxter equation whose semi-classical limit was the original super dynamical $r$-matrix. In particular we have not explicitly constructed algebraic structures which should be the corresponding dynamical quantum groups associated to the resulting $R$-matrices. However, while working in the quantum picture, we have proposed and used a particular algebraic condition which we called {\em the super Hecke condition} (cf.\ Subsection \ref{SS:SuperHecke}). Finding the correct super Hecke condition is important because the Hecke condition in the non-graded case turns out to be the right pre-condition for a meaningful description of dynamical quantum groups in the language of Hopf algebroids (cf.\ \cite{EV2}). 

Studying the proof of our main classification result for super dynamical $R$-matrices (Theorem \ref{T:ClassifyR1}), one can see that the building blocks fall into their right places when one defines the super Hecke condition as we do. In this framework, the super dynamical $R$-matrices with equal Hecke parameters correspond precisely to the zero weight super dynamical $r$-matrices with zero coupling constant. This is exactly analogous to the non-graded picture in \cite{EV2}. This observation may offer some support for our particular definition of the super Hecke condition. 

The construction of the actual algebraic structures that correspond to the super dynamical $R$-matrices we study in Section \ref{S:QuantPict} involves the difficult problem of determining what the appropriate super analogue to dynamical quantum groups should be. This is beyond the scope of this note, but we believe that our work here will shed some light to it by contributing some evidence for the right way to superize the Hecke condition. We intend to address this issue in depth in our followup work. For various possible approaches to the theory of super dynamical quantum groups and some preliminary results, see \cite{Kar4, Kar5}.

\end{document}